\documentclass[12pt]{article}
\usepackage{amssymb,amsfonts,amsmath, psfrag}
\usepackage{graphicx, colordvi}

\newtheorem{theo}{Theorem}

\newtheorem{coro}[theo]{Corollary}

\makeatletter \@addtoreset{equation}{section}
\@addtoreset{theo}{section} \makeatother

\def\qed{\hfill \rule{4pt}{7pt}}
\def\pf{\noindent {\it Proof.} }

\def\CT{{\rm CT\,}}
\def\ep{{\ \equiv_p\ }}

\begin{document}

\begin{center}
{\large \bf
Automated Discovery and Proof of Congruence Theorems for
Partial Sums of Combinatorial Sequences
}

William Y.C. Chen$^{1}$, Qing-Hu Hou$^{2}$, and Doron Zeilberger$^{3}$

$^{1,2}$ Center for Applied Mathematics, Tianjin University\\
Tianjin 300072, P.R. China \\
chenyc@tju.edu.cn, qh\_hou@tju.edu.cn

$^{3}$ Department of Mathematics, Rutgers University (New Brunswick)\\
Piscataway, NJ 08854, USA \\
zeilberg@math.rutgers.edu

\end{center}

\begin{abstract}
Many combinatorial sequences (for example, the Catalan and Motzkin numbers) may be expressed as
the constant term of $P(x)^k Q(x)$,
for some Laurent polynomials $P(x)$ and $Q(x)$ in the variable $x$
with integer coefficients.
Denoting such a sequence by  $a_k$,
we obtain a general formula that determines the congruence class, modulo $p$,
of the indefinite sum  $\sum_{k=0}^{rp -1} a_k$, for {\it any} prime $p$, and any positive integer $r$,
as a linear combination of sequences that satisfy linear recurrence (alias difference) equations with constant coefficients.
This enables us (or rather, our computers) to
automatically discover and prove congruence theorems for such partial sums.
Moreover, we show that in many cases, the set of the residues is finite, regardless of the prime $p$.
\end{abstract}

\section{Introduction}

Let $\{a_k\}$ be a sequence of integers, and $r$ be a positive integer.
We focus on the congruences of the partial sum $\sum_{k=0}^{rp -1} a_k$ modulo a general prime $p$.
When $a_k$ is a hypergeometric term and $r=1$ we get a truncated hypergeometric series, which is closely related to the Gaussian hypergeometric series introduced by Greene \cite{Gre87}. An interesting example is the congruence of the Ap\'ery number \cite{Ahl00,Beu87}
\[
A \left( \frac{p-1}{2} \right) \equiv  \sum_{k=0}^{p-1} {2k \choose k}^4 2^{-8k} \equiv \gamma(p) \pmod{p},
\]
where
\[
A(n) = \sum_{j=0}^n {n+j \choose j}^2 {n \choose j}^2,
\]
and
\[
\sum_{n=1}^\infty \gamma(n)q^n = q \prod_{n=1}^\infty (1-q^{2n})^4 (1-q^{4n})^4.
\]

These congruences are usually obtained case by case and the proofs are complicated. For example, Pan and Sun \cite{Pan06} used a sophisticated combinatorial identity to deduce that
\[
\sum_{k=0}^{p-1} {2k \choose k+d} \equiv \left( \frac{p-d}{3} \right) \pmod{p}, \quad d=0,1,\ldots, p,
\]
where $\left( \frac{\cdot}{\cdot} \right)$ is the Legendre symbol.
We propose an automated method to discover and prove such congruences for a large family of combinatorial sequences $\{a_k\}$. More precisely,
we assume that $a_k$ is the constant term of $P(x)^k Q(x)$ where $P(x)$ and $Q(x)$
are two Laurent polynomials in the (single) variable $x$ with integer coefficients.
Rowland and Zeilberger \cite{Row14} discovered an algorithm to
automatically generate {\it automata} for determining the congruences, modulo a prime $p$,
of combinatorial sequences (not the partial sums) but for {\it specific} primes $p$ (one at a time).

Throughout the paper, $p$ always denotes a prime number.
We write $a \ep b$ if  $a$ is congruent to $b$ modulo $p$.
For a Laurent series $f(x)=\sum_{k \ge k_0} a_k x^k$, we use $\CT f(x)$ to denote the coefficient of the free term, $x^0$.
The set of integers, rational numbers and complex numbers are denoted by $\mathbb{Z}, \mathbb{Q}$ and $\mathbb{C}$, respectively.
 The finite field with $p$ elements is denoted by $\mathbb{F}_p$.

\section{Evaluation}
In this section, we show that the above-mentioned
partial sums are linear combinations of $C$-finite sequences, i.e., integer sequences
that satisfy a linear recurrence equation with {\it constant} coefficients (like $2^n$ and the Fibonacci numbers, to name two examples).
This would enable us (and our computers) to discover and prove practically infinitely-many
theorems about the congruences of such  partial sums modulo an {\it arbitrary} (symbolic!) prime $p$.

We have the following formula for the congruences of the partial sums.
\begin{theo}\label{th-main}
Let $P(x),Q(x)$ be two Laurent polynomials in $x$ with integer coefficients and
\[
a_k := \CT P(x)^kQ(x).
\]
Let $-m$ and $-n$ be the
lowest degrees of $P(x)$ and $Q(x)/(P(x)-1)$, respectively. Then for any positive integer $r$, and any prime $p > n$, we have
\begin{equation}\label{r-case}
\sum_{k=0}^{r p-1} a_k \ep \sum_{j=0}^{rm} c_j S_{(rm-j)p} \ ,
\end{equation}
where $c_j$ is the coefficient of $x^{-rm+j}$ in $P(x)^r - 1$ and $S_k$ is the coefficient of $x^k$ in the Laurent expansion of
\[
\frac{Q(x)}{P(x)-1}.
\]
\end{theo}

\pf
Noting that $\CT$ is a linear operator, we have
\[
\sum_{k=0}^{rp-1} a_k = \CT \sum_{k=0}^{rp-1}  P(x)^k Q(x) =\CT \frac{(P(x)^{rp} - 1)Q(x)}{P(x)-1}.
\]
Since the coefficients of  $P(x)$ are integers, we have, (by the ``{\it Freshman's Dream Identity}'' , $(a+b)^p \ep a^p + b^p$),
$P(x)^p \ep P(x^p)$ and hence
\[
\sum_{k=0}^{rp-1} a_k \ep \CT \frac{Q(x)(P(x^p)^r - 1)}{P(x)-1}.
\]
By the definition of $m$ and $c_j$, we see that
\[
P(x)^r - 1 = \sum_{j=0}^N c_j x^{-rm+j}.
\]
for some integer $N$. If $j>rm$, we have
\[
(-rm+j)p - n > n(-rm+j)-n \ge  0,
\]
which implies that
\[
\CT \frac{Q(x)x^{(-rm+j)p}}{P(x)-1}  = 0.
\]
Hence
\[
\CT  \frac{Q(x)(P(x^p)^r - 1)}{P(x)-1} = \sum_{j=0}^{rm} c_j \CT \frac{Q(x) x^{(-rm+j)p}}{P(x)-1}
= \sum_{j=0}^{rm} c_j S_{(rm-j)p}.
\]
This completes the proof. \qed

This theorem is implemented in the {\tt Maple} package {\tt CTcong.txt} available from the webpage
\begin{center}
{\tt http://www.math.rutgers.edu/\~{}zeilberg/mamarim/mamarimhtml/ctcong.html}
\end{center}
where the user can also find  sample input and output files.

The {\tt Maple} command-line is
\begin{center}
{\tt TheoG(P, Q, x, p, C, r)},
\end{center}
where $P,Q$ are two Laurent polynomials, with integer coefficients, in the variable $x$, $p$ is the symbol standing for the prime,
$C$ is the name for the sequence of coefficients of $Q(x)/(P(x)-1)$, while $r$ is as in Equation~\eqref{r-case}.
For example, typing (in a Maple session, after {\tt reading} our Maple package {\tt CTcong.txt})
\begin{center}
{\tt TheoG(1/x+2+x, x\^\,d, x, p, C, 1);}
\end{center}
immediately outputs
\begin{coro} \label{c1}
Let $A(i)$ be the constant term of the Laurent polynomial
\[
x^d \left( \frac{1}{x} + 2 + x \right)^i,
\]
and for any prime $p$, let
\[
B(p) = \sum_{i=0}^{p-1} A(i).
\]
Then
\[
B(p) \ep C(p),
\]
where $C(n)$ is the $C$-finite sequence defined in terms of the
generating function
\[
\sum_{i=0}^\infty C(i) x^i = \frac{x^{d+1}}{x^2+x+1}.
\]
\end{coro}
Noting that
\[
\CT x^d \left( \frac{1}{x} + 2 + x \right)^i = {2i \choose i-d} = {2i \choose i+d},
\]
and
\[
\frac{x^{d+1}}{x^2+x+1} = x^{d+1} (1+x^3+x^6+\cdots - x - x^4 - x^7 -\cdots),
\]
Corollary~\ref{c1} is equivalent to the congruence relation given by Pan and Sun
\[
\sum_{k=0}^{p-1} {2k \choose k+d} \ep \left( \frac{p-d}{3} \right).
\]

[Of course this case, and many other ones, for {\it small} $r$, are easily humanly-generated.]

Using this approach, we found many congruences, including the congruences for the sums of generalized central trinomial coefficients 
that were considered by Sun in \cite{Sun14}.

When $Q(x)/(P(x)-1)$ is a rational function such that every root of the denominator is a root of unity,
the coefficient of $x^k$ in $Q(x)/(P(x)-1)$
can be expressed as a {\it  quasi-polynomial} in $k$.
We can search for this quasi-polynomial by the method of undetermined coefficients and thus derive theorems presenting {\it explicit}
forms for the congruences.
This is implemented by the procedure {\tt TheoQP} in our Maple package {\tt CTcong.txt}. The  command-line is
\begin{center}
{\tt TheoQP(P, Q, x, p, r, d)} \  ,
\end{center}
where $P,Q$ are the two Laurent polynomials in $x$, $p$ is the symbol standing for
the prime, $r$ is as above, and  $d$ is the expected degree of the searched quasi-polynomial.
(In practice, you start, optimistically, with $d=0$, and if it fails, you keep increasing $d$ to $1$, then $2$,
until you either find something, or give up.)

For example,  typing
\begin{center}
{\tt TheoQP(1/x+2+x, 1, x, p, 2, 0);}
\end{center}
yields
\begin{coro}
\[
\sum_{k=0}^{2p-1} {2k \choose k} \ep
\begin{cases}
  3, & \mbox{if $p \equiv 1 \pmod{3}$}, \\
  -3, & \mbox{if $p \equiv 2 \pmod{3}$}.
\end{cases}
\]
\end{coro}

For more examples, we refer to the above-mentioned webpage
\begin{center}
{\tt http://www.math.rutgers.edu/\~{}zeilberg/mamarim/mamarimhtml/ctcong.html}
\end{center}

\section{Reduction}
In this section, we consider a further reduction of the coefficients $S_p, S_{2p}, \ldots$ in
Equation~\eqref{r-case}. We find that in some cases, the set $\{S_p \mod{p} \}$ of residues is a finite subset of $\mathbb{Q}$ when $p$ runs over all primes.

First, let us consider the coefficients $S_k$ given by
\[
 \frac{ux + v}{a + bx + cx^2} = \sum_{k=0}^{\infty} S_k x^k,
 \]
where $u,v,a,b,c$ are integers, $a \not= 0$ and $a+bx+cx^2$ is irreducible over $\mathbb{Q}$. Let $\Delta = b^2-4ac$ be the discriminate of $ax^2+bx+c$. Since $a+bx+cx^2$ is irreducible, $\Delta \not= 0$, and hence $\Delta \not\equiv_p 0$ except for finitely many primes $p$.

If $\Delta$ is a square in the finite field $\mathbb{F}_p$, then $a+bx+cx^2$ is reducible in $\mathbb{F}_p$ so that
\[
\frac{ux + v}{a + bx + cx^2} \ep \frac{A}{1-\alpha x} + \frac{B}{1-\beta x},
\]
for some $A,B,\alpha,\beta \in \mathbb{F}_p$. We thus find that
\[
S_{rp}= A \alpha^{rp} + B \beta^{rp} \ep A \alpha^r + B \beta^r = S_r.
\]

If $\Delta$ is not a square in $\mathbb{F}_p$, then $a+bx+cx^2$ is irreducible in $\mathbb{F}_p$.
Let us consider the extension field $\mathbb{F}_p(\alpha)$ with $a \alpha^2 + b \alpha + c = 0$ and
$\alpha \in \mathbb{C}$. Let $\beta \in \mathbb{C}$ be another root of the equation $ax^2+b x + c = 0$.
By the property of the Frobenius automorphism \cite{Wan03}, it follows that in the extension field $\mathbb{F}_p(\alpha)$,
\[
\alpha^p = \beta, \quad \beta^p = \alpha.
\]
Hence in the field $\mathbb{F}_p(\alpha)$, we have
\[
S_{rp} = A \alpha^{rp} + B \beta^{rp} = A \beta^r + B \alpha^r = \frac{c^r}{a^r} \left( A \alpha^{-r} + B \beta^{-r} \right)  = \frac{c^r}{a^r} S_{-r},
\]
where $S_{-r}$ is determined by the initial values $S_0, S_1$ and the recurrence relation
\[
a S_n + b S_{n-1} + c S_{n-2} = 0, \quad n \in \mathbb{Z}.
\]
Since $S_{rp}$ and $S_{-r}$ are both rational numbers, we obtain that $S_{rp} \ep S_{-r}$.

In general, let $q(x)$ be an irreducible polynomial in $\mathbb{Z}[x]$ of degree $d$ with non-zero constant term and $\alpha_1,\ldots, \alpha_d$ be the $d$ roots of $x^d q(1/x)$ in $\mathbb{C}$. If the splitting field $\mathbb{Q}(\alpha_1,\ldots, \alpha_d)$ equals $\mathbb{Q}(\alpha_j)$ for some $1 \le j \le d$, we say that $q(x)$ is {\it simple}. Clearly, every irreducible polynomial of degree $2$ is simple.

We have the following finiteness theorem regarding the congruences.
\begin{theo}\label{finite}
Let $P(x),Q(x)$ be two Laurent polynomials in $x$ with integer coefficients and
\[
a_k=\CT P(x)^kQ(x).
\]
Suppose that each irreducible factor $q(x) \not= x$ of the denominator of $Q(x)/(P(x)-1)$ is simple.
Then there exists a finite subset $A$ of $\mathbb{Q}$ such that for any prime $p$,
\[
\sum_{k=0}^{r p-1} a_k \ep a,
\]
for some $a \in A$.
\end{theo}
\pf
By Theorem~\ref{th-main}, for a sufficiently large prime $p$,
$\sum_{k=0}^{r p-1} a_k$ modulo $p$ is a linear combination of $S_0, S_p, S_{2p}, \ldots$,
where $S_k$ is the coefficient of $x^k$ in the series
\[
R(x) = \frac{Q(x)}{P(x)-1}.
\]
To prove the theorem, it suffices to show that for a fixed integer $n$, the set
\[
\bigcup_p\  \{S_{np} \mod{p}\}
\]
of residues is finite when $p$ runs over all primes.

Consider the partial fraction decomposition of $R(x)$ over $\mathbb{Q}$
\[
R(x) =  g(x) + \sum_{i=1}^m \frac{h_i(x)}{q_i(x)^{\ell_i}},
\]
where $g(x)$ is a Laurent polynomial over $\mathbb{Q}$ and for each $i=1,\ldots,m$, $q_i(x) \in \mathbb{Z}[x]$ is irreducible and $h_i(x) \in \mathbb{Z}$ is a polynomial with ${\rm deg\,} h_i(x)<{\rm deg\,}  q_i(x)$.
In order to show the finiteness of the set $\cup_p \{S_{n p} \mod{p}\}$,
it suffices to show that the residues of the coefficients of each summand form a finite set.

Let $h(x)/q(x)^\ell$ be one summand and
\[
\sum_{k=0}^\infty s_k x^k = \frac{h(x)}{q(x)^{\ell}}.
\]
Let $\tilde{q}(x)=x^d q(1/x)$ where $d = {\rm deg\,} q(x)$ and let $\alpha_1,\ldots,\alpha_d$ be the roots of $\tilde{q}(x)$. By assumption, we have $\mathbb{Q}(\alpha_1,\ldots,\alpha_d) = \mathbb{Q}(\alpha)$ with $\alpha \in \{\alpha_1,\ldots,\alpha_d \}$.
Denote the splitting field $\mathbb{Q}(\alpha)$ by $K$. Since $q(x)$ is irreducible, we have
\[
K = \left\{ \frac{a_0 + a_1 \alpha + \cdots + a_{d-1} \alpha^{d-1}}{b} \colon a_0,\ldots,a_{d-1},\, b \in \mathbb{Z} \right\}.
\]
Let $p$ be a prime larger than the maximal factor of the leading coefficient of $\tilde{q}(x)$.
Then
\[
K_p := \left\{ \frac{a_0 + a_1 \alpha + \cdots + a_{d-1} \alpha^{d-1}}{b} \colon p \nmid b \right\} \subset K
\]
form a subring of $K$. There is a natural ring homomorphism $\tau \colon K_p \to \mathbb{F}[x]/ \langle \tilde{q}(x) \rangle$ given by
\[
\tau \left( \frac{a_0 + a_1 \alpha + \cdots a_{d-1} \alpha^{d-1}}{b} \right) =  a_0 b^{-1} + a_1 b^{-1} x + \cdots + a_{d-1} b^{-1} x^{d-1}.
\]
Clearly, the kernel of the map $\tau$ is $p K_p$.

It is well-known that the coefficients $s_k$ can be expressed as
\[
s_k = \sum_{i=1}^d f_i(k) \alpha_i^k,
\]
where $f_i(k)$ is a polynomial over $\mathbb{Q}(\alpha)$ of degree less than $\ell$.
It is easy to see that
\[
\tau(f_i(np)) = \tau(c_i),
\]
where $c_i$ is the constant term of $f_i(x)$. Since $\tilde{q}(x) \in \mathbb{Z}[x]$, for each $i=1,\ldots,d$,
we have
\[
\tilde{q} (\tau(\alpha_i^{p})) = \tau\big(\tilde{q} (\alpha_i^{p}) \big) = \tau\big( (\tilde{q} (\alpha_i))^p \big) = 0.
\]
Hence $\tau(\alpha_i^{p}) = \tau(\alpha_j)$ for some $1 \le j \le d$. Let $\sigma$ be the map given by $\tau(\alpha_i^{p}) = \tau(\alpha_{\sigma(i)})$. We thus have
\[
\tau(s_{np}) = \sum_{i=1}^d \tau(c_i) \tau(\alpha_{\sigma(i)}^n) = \tau \left( \sum_{i=1}^d c_i \alpha_{\sigma(i)}^n \right).
\]
Let
\[
\sum_{i=1}^d c_i \alpha_{\sigma(i)}^n = r_0 + r_1 \alpha + \cdots + r_{d-1} \alpha^{d-1}.
\]
Since $\tau(s_{np}) \in \mathbb{Q}$, we have
\[
\tau (s_{np}) = \tau(r_0),
\]
and hence $s_{np} \ep r_0$.
Noting that there are only finitely many choices for $\sigma$, hence the set $\cup_p \{s_{np} \mod{p} \}$ is finite.
\qed

\noindent {\it Example.} Suppose that
\[
P(x)=\frac{x^3-2x+1}{x}, \quad \mbox{and} \quad Q(x)=1.
\]
We have
\[
\frac{Q(x)}{P(x)-1} = \frac{x}{x^3-3x+1}.
\]
Let
\[
\alpha=-0.532 \ldots, \quad \beta=0.6527 \ldots, \quad \gamma =  2.879 \ldots
\]
be the three roots of $x^3-3x^2+1$. Using the approximate values of the three roots, we may use the LLL algorithm \cite{Len82} to find integral relations among $\beta,\gamma$ and powers of $\alpha$. Using Maple, we get two possible relations
\begin{equation}\label{eq-abc}
\beta = 2 + 2\alpha - \alpha^2, \quad \gamma = 1 -3\alpha +\alpha^2.
\end{equation}
It is easy to verify that
\[
(2 + 2\alpha - \alpha^2)^3 - 3(2 + 2\alpha - \alpha^2)^2 +1 = 0,
\]
and
\[
(1 -3\alpha +\alpha^2)^3 - 3(1 -3\alpha +\alpha^2)^2 + 1 = 0.
\]
which means that $2 + 2\alpha - \alpha^2$ and $1 -3\alpha +\alpha^2$ are roots of $x^3-3x^2+1$. So we claim the relations in \eqref{eq-abc}.
Therefore, $\mathbb{Q}(\alpha,\beta,\gamma) = \mathbb{Q}(\alpha)$ and hence $x^3-3x+1$ is simple. By Theorem~\ref{finite}, the set
\[
\left\{ \sum_{k=0}^{2p-1} \CT P(x)^k \mod{p} \right\}
\]
is finite. In fact, when $p>3$, the only possibilities are $-1$ and $2$.

We conclude with an example where the denominator is not simple.

\noindent {\it Example}. Let
\[
P(x) = -2x^2 + 1 + \frac{1}{x}, \quad Q(x) =1.
\]
Then
\[
\frac{Q(x)}{P(x)-1} = \frac{x}{1-2x^3} = x + 2x^4 + 2^2 x^7 + \cdots.
\]
Hence for $p \equiv 1 \pmod{3}$, we have
\[
\sum_{k=0}^{p-1} \CT P(x)^k \ep 2^{\frac{p-1}{3}} \ep 2^{-\frac{1}{3}}.
\]
It seems that the set $\{ 2^{-\frac{1}{3}} \mod p\}$ is not finite.

\vskip 15pt
{\bf Acknowledgements.} We wish to thank Zhi-Wei Sun for valuable suggestions. This work was supported by the 973 Project, the PCSIRT Project of the Ministry of Education, and the National Science Foundation of China.

\bibliography{mod}
\bibliographystyle{plain}

\end{document}